\documentclass[12pt]{article}
\usepackage{amssymb}
\usepackage{amsthm}
\usepackage{enumerate}
\usepackage{amssymb}
\begin{document}
\newcommand{\R}{{\bf R}}
\newcommand{\K}{{\bf K}}
\newcommand{\C}{{\bf C}}
\newcommand{\sequ}[1]{(#1_{n})_{n=1,2,\dots}}
\newcommand{\seq}[1]{(#1_{n})_{n=0,1,\dots}}
\newtheorem{thm}{Theorem}[section]
\newtheorem{defn}[thm]{Definition}
\newtheorem{cor}[thm]{Corollary}
\newcommand{\rtref}[1]{{\rm \ref{#1}}}
\newcommand{\rmref}[1]{{\rm (\ref{#1})}}
\newcommand{\rmcite}[1]{{\rm \cite{#1}}}
\newcommand{\rmciter}[2]{{\rm \cite[#1]{#2}}}
\newtheorem{lem}[thm]{Lemma}
\newtheorem{prop}[thm]{Proposition}
\newtheorem{rem}[thm]{Remark}

\baselineskip=1.5\baselineskip




\title {On the Decreasing Failure Rate property for general counting process. Results based on conditional interarrival times}

\author {F. G. Bad\'{\i}a, C.Sang\"{u}esa
 \\ {\small
Departamento de M\'etodos Estad\'{\i}sticos and IUMA. Universidad de
Zaragoza. }
 \\ {\small Zaragoza. Spain. } \\ {\small e-mail: gbadia@unizar.es,csangues@unizar.es} }

\date{}

\begin{titlepage}
\setcounter{page}{1} \maketitle

\bigskip \bigskip
\begin{abstract}
In the present paper we consider general counting processes stopped at a random time $T$, independent of the process.  Provided that $T$ has the decreasing failure rate (DFR) property, we give sufficient conditions on the arrival times so that the number of events occurring
before $T$ preserves the DFR property of $T$. These conditions involve the study of the conditional interarrival times.  As a main application, we prove the DFR property in a context of maintenance models in reliability, by the consideration of Kijima type I virtual age models under quite general assumptions.
\end{abstract}



\end{titlepage}

\section{Introduction}
Decreasing failure rate (DFR) is a property describing a system which improves with age.
For instance, the lifetime of a mechanism whose distribution is a mixture of exponential random variables has this property
(see Proschan \cite{prtheo}). In fact, it is well known that mixtures of DFR distribution have the DFR property (see also \cite{guseho} or \cite{fieswh}, for instance, for different contexts in which the DFR property can arise).  Motivated by the fact that the DFR property is easily preserved by mixtures, and with the aim of completing the previous results obtained in \cite{basath}, we will study the discrete DFR property for $N(T)$, in which $ \{ N(t): t \geq 0 \} $ is a counting process and $T$ is a random time independent of this process.  The discrete DFR property is a logconvexity condition for the survival function (see Esary {\it et al.}
\cite{esmash} or Grandell \cite[Ch. 7]{grmixe}), that is,
 \begin{eqnarray} P(N(T)\geq n+1)^{2}&\leq&
P(N(T)\geq n)P(N(T)\geq n+2) \quad n=0,1,\dots . \label{disdfr}
\end{eqnarray}

Note that in this case we consider a mixture of random variables which, in general don't have the DFR property (for instance, in a Poisson process, each $N(t)$ is discrete increasing failure rate).

The interest in the study of $N(T)$ comes back to \cite{corene}, p. 42 in which several examples were given. Natural applications arise also in queuing theory, when studying the stationary number of waiting customers (see \cite{assamp}, for instance).  Recent applications have been found in health sciences (see \cite{soesti} ), in which the counting process describes the number of tests for a disease (HIVS for instance) in a person at risk until the random time $T$ at which this person is infected.

Preservation of different ageing properties for a classical renewal process (that is, with independent, and identically distributed interarrival times) have been extensively studied (see \cite{esmash,grmixe,BS,basapr}, for instance).  Recently, attention is focused in generalizing these properties for counting processes exhibiting dependence between the interarrival periods (see \cite{roshon,german}). These models are of interest in reliability, for instance, when considering a process which models the successive repairments of a mechanism which deteriorates with age (see \cite{kisuau}). In the study of these models, it is sometimes usual to consider the simplifying hypothesis that the interarrival epochs are discrete or absolutely continuous random vectors and one of the aims of this paper is to present our results in a general background easily adapted to usual examples.

As mentioned before, we complete the DFR preservation results obtained in \cite{basath}, going deeper in the study of the dependence properties between the interarrival times which can give place to this property,  with special attention to get preservation results for general counting models of interest in reliability. Our main result in this paper is Proposition \ref{t2*}, which shows the preservation of the DFR property under a (strong) decreasing behaviour of the interarrival epochs. The result is based on Propositions \ref{prdfr2} and \ref{prcon}.  The last one provides a way to check (\ref{disdfr}) for $n=1$, when the second renewal epoch satisfy our decreasing condition.  Proposition \ref{prdfr2} provides a way to check (\ref{disdfr}) for a general $n$ by means of the conditional distribution of $(X_{n+1},X_{n+2})$ given the past.  Our result differs from the ones given in \cite{basath}, in the sense that the last ones need, along with (different) decreasing conditions for the renewal epochs, some additional conditions. As an application we will show in Section 3 the DFR property for Kijima type I models (see \cite{kisuau}) under more general assumptions than in \cite[Section 5]{basath}.  Finally, as an immediate application of Proposition  \ref{prcon} we will give another preservation result which assumes the condition used in \cite[Proposition 3.1.(a)]{basath}) for the first interarrival epochs, together with the same condition for the the next two renewal epochs knowing the past history.

The proofs of our results are made without considering the absolutely continuity assumption on the interarrival epochs before mentioned, so that in Section \ref{prelim} we introduce some preliminar background in order to specify the conditional distribution for two inter renewal times given the past history, making use of conditional probability kernels. Also in Section \ref{prelim}, we give some preliminary definitions which are going to be used along the paper.

\section{Preliminaries} \label{prelim}
Let $ \{ N(t): t \geq 0 \} $ be a  general counting process whose arrival times are denoted by $S_n$, $ n=0,1, \ldots  \  (S_0=0)$ and the interarrival epochs are denoted by $X_n$, $n=1, 2, \ldots $, that is:
 \[
 S_n = \sum_{i=1}^n X_i, \quad n=1, 2, \ldots
 \]
 The only assumption about the $X_i$ is that they are non negative random variables. Recall that the counting process is defined through the renewal epochs by the following expression
 \[
 N(t) = \max \{n: S_n \leq t \},  \quad t \geq 0.
 \]

Our paper deals with the preservation of the decreasing failure rate property. This concept, for discrete random variables, was reminded in the Introduction.  The concept for 'general' random variables is recalled below (see Barlow and Proschan
\cite{baprst}) along with another one reliability class (New Worse than Used), which is weaker than the DFR property, and will be used throughout the paper.
 \begin{defn} Let $X$ be a nonnegative random variable with $G$ and \mbox{$\overline
 G:=1-G$} the corresponding distribution and reliability
 functions.  $X$  (or $G$) is said to be:
\begin{enumerate} \item \emph{Decreasing failure  rate}  if

$\overline G (z+t)/\overline G(t)$ is increasing in $t$, for all $z\geq 0$;

\item \emph{New worse than used}
(NWU)  if

$\overline{G}(z)\overline{G}(t)\leq\overline{G}(z+t)$, for all
$z,t\geq 0$;

\end{enumerate}\label{dfgrel}\end{defn}

Note that if the monotonicity and sense of the inequalities in Definition \ref{dfgrel} are reversed, we obtain the dual concepts of \emph{increasing failure rate} (IFR) and \emph{new better than used}.

Our results make use of comparisons of random variables in the usual stochastic order, which we recall in the following.
\begin{defn} Let $X$ and $Y$ be two random variables with distribution functions $F_{X}$ and $F_{Y}$, respectively. $X$ is said to be \emph{smaller
than $Y$ in the usual stochastic order} (denoted $X\leq_{ST} Y$) if $F_{X}(t)\geq F_{Y}(t)$ for all real  $t$ or, equivalently, if $\overline{F}_{X}(t)\leq\overline{F}_{Y}(t)$ for all real  $t$.\label{dfstor}
\end{defn}

In our results concerning the DFR preservation property in a renewal process, we will make use of the conditional distributions of the inter renewal epochs. More specifically, for a given natural number $m$, we will consider the conditional distribution of $(X_{n+1},\dots,X_{n+m})$ given that $(X_{1},\dots X_{n})$ is known. In order to deal in a general setting, we will use conditional probability kernels. Consider two Borel sets $S\subseteq \R^{n}$ and $T\subseteq \R^{m}$.  Denote by ${\mathcal B}(T)$ the family of Borel sets on $T$.  Recall that a probability kernel $\mu$ from $S$ to $T$ is a mapping \begin{eqnarray*}\mu: S\times {\mathcal B}(T)&\rightarrow & [0,1]\\(\mathbf{x},B)&\rightsquigarrow& \mu^{\mathbf{x}}(B)\end{eqnarray*}  verifying the following properties:
\begin{enumerate}
\item $\mu^{\mathbf{\cdotp}}(B)$ is a measurable function for $B\in {\mathcal B}(T)$ fixed.
\item $\mu^{\mathbf{x}}(\cdotp)$ is a probability distribution on $(T,{\mathcal B}(T))$, for  $\mathbf{x}\in S$ fixed.
\end{enumerate}
Let  $(X_{1},\dots,X_{n},X_{n+1},\dots X_{n+m})$ be a nonnegative random vector defined on a given probability space.   Taking into account \cite[Thm. 6.3, p.107]{kafoun}), the existence of a probability kernel $\mu$ from $\R^{n}_{+}$ to $\R^{m}_{+}$ such that for each Borel set $B\subseteq \R_{+}^{m}$, $\mu^{(X_{1},\dots,X_{n})}(B)$ is a version of $P((X_{n+1},\dots,X_{n+m})\in B\vert(X_{1},\dots,X_{n}))$ can be always guaranteed.  Moreover, if we denote by $F_{n}$ the distribution function of $(X_{1},\dots,X_{n})$, we have for a nonnegative measurable  function $f:\R_{+}^{n}\rightarrow \R_{+}$ that (cf. \cite[p.108]{kafoun})

\[E[f(X_{1},\dots,X_{n+m})]=\int_{\mathbf{z}\in \R_{+}^{n+m}}f(\mathbf{z})dF_{n+m}(\mathbf{z})=
\int_{\mathbf{x}\in\R_{+}^{n}}dF_{n}(\mathbf{x})\int_{\mathbf{y}\in \R_{+}^{m}} f(\mathbf{x},\mathbf{y})
d\mu^{\mathbf{x}}(\mathbf{y}) \]
and obviously, if we take a set $N_{n}\in {\mathcal B}(\R_{+}^{n})$ such that $P((X_{1},\dots,X_{n})\in N_{n})=1$, we have
\begin{equation}E[f(X_{1},\dots,X_{n+m})]=\int_{\mathbf{x}\in N_{n}}dF_{n}(\mathbf{x})\int_{\mathbf{y}\in \R_{+}^{m}} f(\mathbf{x},\mathbf{y})
d\mu^{\mathbf{x}}(\mathbf{y}). \label{prodint}\end{equation}
In fact, (\ref{prodint}) is the property we will need in our proofs.  Observe that (\ref{prodint}) is satisfied for all nonnegative $f$ if an only if it is satisfied for $f(\mathbf{x},\mathbf{y})=1_{A\times B}(\mathbf{x},\mathbf{y})$, for Borel sets $A\subseteq N_{n},\ B\subseteq \R_{+}^{m}$, in which $1_{A\times B}(\cdotp)$ denotes the indicator function on the set $A\times B$ (see \cite[Theorem 2.6.4, p.105]{asdopr}). This motivates the following:
\begin{defn} Consider a random vector $(X_{1},\dots,X_{n},X_{n+1},\dots X_{n+m})$ of nonnegative random variables. Consider $N_{n}\in {\mathcal B}(\R_{+}^{n})$ such that $P((X_{1},\dots,X_{n})\in N_{n})=1$. Let $\mu$ be a probability kernel from $N_{n}$ to $\R^{m}_{+}$ verifying for each Borel sets $A\subseteq N_{n},\ B\subseteq\R_{+}^{m}$,
\begin{equation}\int_{\{\mathbf{x}\in A,\mathbf{y}\in B\}}dF_{n+m}(\mathbf{x},\mathbf{y})=\int_{A}dF_{n}(\mathbf{x})\mu^{\mathbf{x}}(B).\label{promed}\end{equation}

The family $\{\mu^{\mathbf{x}}(\cdotp),\ \mathbf{x}\in N_{n}\}$
will be said a $N_{n}$ regular conditional distribution of $(X_{n+1},\dots,X_{n+m})$ given $(X_{1},\dots,X_{n})$. A family of m-dimensional random vectors $\{(Z_{n+1}^{ \mathbf{x}},\dots, Z_{n+m}^{ \mathbf{x}}),\ \mathbf{x}\in N_{n} \}$,  such that for each $\mathbf{x}\in N_{n}$,  the distribution of $(Z_{n+1}^{ \mathbf{x}},\dots, Z_{n+m}^{ \mathbf{x}})$ is given by $\mu^{\mathbf{x}}(\cdotp)$,
will be said a $N_{n}$-distributional version of
$(X_{n+1},\dots,X_{n+m})$ given $(X_{1},\dots,X_{n})$.
\label{defcon}\end{defn}

\begin{rem} The following remarks are in order:
\begin{enumerate}\item If in the previous definition the distribution of $(X_{1},\dots,X_{n+m})$ is discrete (resp. absolutely continuous), the set $N_{n}$ can be taken as the set of points $\mathbf{x}\in \R_{+}^{n}$ such that $(X_{1},\dots,X_{n})=\mathbf{x}$ has strictly positive probability (resp. strictly positive density function), and the conditional distributions, defined in the usual way, satisfy (\ref{promed}). Another important example is when $(X_{1},\dots,X_{n})$ and  $(X_{n+1},\dots,X_{n+m})$ are independent. In this case we can take $N_{n}=\R^{n}_{+}$ and $\mu^{\mathbf{x}}(\cdotp)$, constant for all $\mathbf{x}$, being the distribution of $(X_{n+1},\dots,X_{n+m})$.
\item The random vectors $(Z_{n+1}^{ \mathbf{x}},\dots, Z_{n+m}^{ \mathbf{x}})$  in Definition \ref{defcon} are introduced for notational convenience (actually, they don't need to be defined on the same probability space, for different values of $\mathbf{x}$).  In fact, the properties we need to use (stochastic monotonicity, for instance) will apply to the vector $(Z_{n+1}^{ \mathbf{x}},\dots, Z_{n+m}^{ \mathbf{x}})$, with $\mathbf{x}$ fixed and have to do with the probability distribution of this vector, not with the paths of the process.
  \end{enumerate} \label{redeco}\end{rem}

\section{DFR property under decreasing conditions of the interarrival times}\label{second}
In this Section we present our first result concerning the DFR preservation property in a renewal process, making use of the conditional distributions of the inter renewal epochs.
  Let $  \{ N(t): t \geq 0 \} $ be a general counting process with inter renewal epochs $\seq{X}$. Let $T$ be a DFR random time independent from the process, whose survival function is denoted by $\overline {F}_T$. We know  that if $T$ and $X_1$ don't have 0 as a common discontinuity point, then (cf. \cite[Lemma 3.2]{basath})
\begin{equation}
P(N(T) \geq n ) = E (\overline {F}_T (S_n)), \quad n =1,2,\dots \label{nocomm}
\end{equation}and inequality (\ref{disdfr}) can be rewritten, for each $n$ as:
  \begin{eqnarray}E^{2} (\overline {F}_T (X_{1}))&\leq& E [\overline {F}_T (S_{2})],\quad n=0,\label{dfrn0}\\
E ^{2}(\overline {F}_T (S_{n+1}))&\leq &E (\overline {F}_T (S_n))E (\overline {F}_T (S_{n+2})),\quad n=1,2,\dots .\label{dfrn}\end{eqnarray}

 \medskip In the next result we give sufficient conditions in order to verify (\ref{dfrn}) when we can find distributional  versions of $(X_{n+1},X_{n+2})$ given $(X_{1},\dots,X_{n})$ satisfying (\ref{dfrn0}).
\begin{prop}
Let $  \{ N(t): t \geq 0 \} $ be a general counting process with inter renewal epochs $\seq{X}$. Let $T$ be a DFR random time independent from the process. Let $n$ be a fixed natural number and assume that there exists a $ N_{n}$ distributional version of of $(X_{n+1},X_{n+2})$ given $(X_{1},\dots,X_{n})$ (say $\{(Z_{n+1}^{\mathbf{x}},Z_{n+2}^{\mathbf{x}}),\quad  \mathbf{x}\in N_{n}\}$) satisfying, for every DFR survival function  $\overline {H}$ with  $\overline {H}(0)=1$, that
\begin{equation}E^{2}[ \overline {H} (Z_{n+1}^{\mathbf{x}})] \leq E[ \overline {H}( Z_{n+1}^\mathbf{x}+Z_{n+2}^{\mathbf{x}})],\quad \mathbf{x}\in N_{n}. \label{concon}\end{equation}
Then,
\begin{equation} E ^{2}(\overline {F}_T (S_{n+1}))\leq E (\overline {F}_T (S_n))E (\overline {F}_T (S_{n+2})). \label{dfrn*}\end{equation}
\label{prdfr2}
\end{prop}

{\bf Proof:}  In order to check inequality (\ref{dfrn*}), consider $\mathbf{x}:=(x_{1},\dots,x_{n})\in N_{n}$ and call $s_{n}:=x_{1}+\dots+x_{n}$.  As $\overline {F}_T $ is a DFR survival function, then the survival function $\overline {H}$ defined as
\begin{equation}\overline {H}(z)=\frac{\overline {F}_T (s_{n}+z)}{\overline {F}_T (s_{n})},\quad z\geq 0\label{auxh}\end{equation}
verifies the DFR property (see \cite[p.118]{maolli})) and, as $\overline {H}(0)=1$, condition (\ref{concon}) holds true. If we call $\mu^{\mathbf{x}}$ to the common distribution of $(Z_{n+1}^{\mathbf{x}}, Z_{n+2}^{\mathbf{x}})$, we can write (\ref{concon}) as
\[\int_{\R_{+}^{2}} \overline {H}(x_{n+1})d\mu^{\mathbf{x}}(x_{n+1},x_{n+2})\leq\left(\int_{\R_{+}^{2}} \overline {H}(x_{n+1}+x_{n+2})d\mu^{\mathbf{x}}(x_{n+1},x_{n+2})\right)^{1/2}\]
which, recalling (\ref{auxh}), turns out to  be
\begin{eqnarray}&&\int_{\R_{+}^{2}} \overline {F}_T(s_{n}+x_{n+1})d\mu^{\mathbf{x}}(x_{n+1},x_{n+2})\nonumber\\&&\leq\overline {F}_T^{1/2}(s_{n})\left(\int_{\R_{+}^{2}} \overline {F}_T(s_{n}+x_{n+1}+x_{n+2})d\mu^{\mathbf{x}}(x_{n+1},x_{n+2})\right)^{1/2},\label{chorizo}\end{eqnarray}
Call
\begin{equation}g(\mathbf{x}):=\int_{\R_{+}^{2}} \overline {F}_T(s_{n}+x_{n+1}+x_{n+2})d\mu^{\mathbf{x}}(x_{n+1},x_{n+2}),\quad \mathbf{x}\in N_{n}. \label{gesene}
\end{equation}
Integrating (\ref{chorizo}) with respect to $F_{n}$, the distribution function of $(X_{1},\dots,X_{n})$, we have

\begin{equation}\kern-20pt\int_{N_{n}}dF_{n}(\mathbf{x})\int_{\R_{+}^{2}} \overline {F}_T(s_{n}+x_{n+1})d\mu^{\mathbf{x}}(x_{n+1},x_{n+2})\leq \int_{N_{n}}dF_{n}(\mathbf{x})\overline {F}_T^{1/2}(s_{n})g^{1/2}(\mathbf{x}).\label{totoro}\end{equation}
Applying Cauchy-Schwartz inequality to the term in the right-hand side of (\ref{totoro}) we obtain

\begin{equation}\kern-20pt \int_{N_{n}}dF_{n}(\mathbf{x})\overline {F}_T^{1/2}(s_{n})g^{1/2}(\mathbf{x})\leq\left(\int_{N_{n}}dF_{n}(\mathbf{x})\overline {F}_T(s_{n})\right)^{1/2}\left(\int_{N_{n}}dF_{n}(\mathbf{x})g(\mathbf{x})\right)^{1/2}.\label{t*aux1}\end{equation}
Thus, by (\ref{totoro}), (\ref{t*aux1}) and recalling (\ref{gesene}) we can write
\begin{eqnarray*}&&\kern-30pt \left(\int_{N_{n}}dF_{n}(\mathbf{x})\int_{\R_{+}^{2}} \overline {F}_T(s_{n}+x_{n+1})d\mu^{\mathbf{x}}(x_{n+1},x_{n+2})\right)^{2}\\&&\kern-30pt\leq
\left(\int_{N_{n}}dF_{n}(\mathbf{x})\overline {F}_T(s_{n})\right)\left(\int_{N_{n}}dF_{n}(\mathbf{x})\int_{\R_{+}^{2}} \overline {F}_T(s_{n}+x_{n+1}+x_{n+2})d\mu^{\mathbf{x}}(x_{n+1},x_{n+2})\right).\end{eqnarray*}
The previous inequality shows (\ref{dfrn*})$\quad \Box$

Two results concerning the preservation of the DFR property were given in \cite{basath}.  The first result (See Section 3 in \cite{basath})  involves the association (see Definition \ref{prasso} in Section 5 below) and stochastic decrease of the interarrival times.  In the second result (see Section 5 in \cite{basath}) a different set of conditions were given. These conditions also involve a decreasing behaviour of the interarrival times (together with some other technical conditions).  We now introduce a different way to check (\ref{dfrn0}) and (\ref{dfrn}), assuming a stronger decreasing condition on the interarrival times than in the two previous results, but with no more additional conditions.  This property will allow us to check the d-DFR condition for generalized renewal processes under more general assumptions than in \cite[Section 5]{basath}.
\begin{prop} Let $(X_1,X_2)$ be a nonnegative random vector.  Assume that $ \overline {F}_T $ is a DFR survival function.  Assume  also, that there exists a $N_{1}$ distributional version of $X_{2}$ given $X_1$ (say $\{Z_{2}^{x_{1}}, \ x_{1}\in N_{1}\})$, such that $Z_{2}^{x_{1}}\leq_{ST}X_{1}$, for all $x_{1}\in N_{1}$. Then,
\[E^2[ \overline {F}_T (X_{1})] \leq E[ \overline {F}_T (X_{1}+X_{2})]\overline {F}_T (0).  \]
\label{prcon}
\end{prop}
{\bf Proof:}  First of all, if $\overline {F}_T$ is is a DFR survival function, then $\overline {F}_T(\cdotp)/\overline {F}_T(0)$ is a DFR survival function and, therefore, NWU (see \cite[p.181]{maolli}).  Thus,
\[
\frac{\overline {F}_T(z+t)}{\overline {F}_T(0)} \geq \frac{\overline {F}_T(z)\overline {F}_T(t)}{(\overline {F}_T(0))^{2}}, \quad z, t \geq 0.
\]
Then, choosing $z=X_1$ and $t=X_2$ and taking expectations
\begin{equation}
\overline{F}_T(0) E(\overline {F}_T (X_1+X_2)) \geq E(\overline {F}_T (X_1)\overline {F}_T (X_2))
\label{laux1}\end{equation}
Now, consider the family $\{Z_{2}^{x_{1}}, \ x_{1}\in N_{1}\}$, and call $\{\mu_{1}^ {x_{1}}, \ x_{1}\in N_{1}\}$ to its associate family of distributions.  Recalling (\ref{prodint}), we can write
\begin{eqnarray}
E(  \overline {F}_T (X_1)\overline {F}_T (X_2))&=&\int_{N_{1}}\overline {F}_T (x_1)dF_{1}(x_{1})\int_{R^{+}}\overline {F}_T (x_{2})d\mu_{1}^ {x_{1}}(x_{2})\nonumber\\&=&\int_{N_{1}}\overline {F}_T (x_1)E(\overline {F}_T (Z_{2}^{x_{1}}))dF_{1}(x_{1}).
\label{laux1b}
\end{eqnarray}
As we have in this case that $Z_{2}^{x_{1}}\leq_{ST}X_{1}$, then
\begin{equation}E(\overline {F}_T (Z_{2}^{x_{1}}))\geq E(\overline {F}_T (X_1)).\label{laux2b}\end{equation}
Thus, we deduce from (\ref{laux1b}) and (\ref{laux2b})
\[E(  \overline {F}_T (X_1)\overline {F}_T (X_2))\geq E^{2}(\overline {F}_T (X_1)).\]
This inequality together with (\ref{laux1}) prove Proposition \ref{prcon}.$\quad \Box$

\medskip Our aim now is to show the preservation of the DFR property making use of Propositions \ref{prdfr2} and \ref{prcon}. This involves to find a distributional version of $(X_{n+1},X_{n+2})$ given $(X_{1},\dots,X_{n})$ satisfying the conditions in Proposition \ref{prcon}.  However, in many situations, unidimensional conditional distributions  are given in a natural way, that is, we know the distribution of $X_{n+1}$ given $(X_{1},\dots,X_{n}),\ n=1,2,\dots$ . In the discrete, or absolutely continuous case, the construction of bidimensional conditionals provided sequentially unidimensional conditionals is well-known. In order to extend this formula, recall that if we have a probability kernel $\mu_{n+1}$ from $N_{n}\in \ {\mathcal B}(R_{+}^{n})$ to $R_{+}$, and we have a probability kernel $\mu_{n+2}$ from $N_{n}\times R_{+}$ to $R_{+}$, we can define the composition kernel $\nu:=\mu_{n+1}\otimes\mu_{n+2}$ from $N_{n}$ to $\R^{2}_{+}$ as (cf. \cite[p. 20]{kafoun}
\begin{equation} \nu_{n}^{\mathbf{x}}(B_{1}\times B_{2})=\int_{B_{1}}d\mu_{n}^{\mathbf{x}}(x_{n+1})\int_{B_{2}}d\mu_{n+1}^{(\mathbf{x},x_{n+1})}(x_{n+2}),\quad \mathbf{x}\in N_{n}, \label{condos}\end{equation}
for any Borel sets $B_{1}\in {\mathcal B}(\R_{+})$ and $B_{2}\in {\mathcal B}(\R_{+})$.  In the next Proposition we show that composition kernels are the natural way to construct bidimensional conditional probability kernels. Its proof, although using standard techniques, is rather long, so that it is postponed to the appendix.
\begin{prop} Let $\seq{X}$ be a family of inter renewal epochs. For $n$ fixed, call $\mathbf{X}=(X_{1},\dots,X_{n})$ and assume that $\{\mu_{n}^{\mathbf{x}},\ \mathbf{x}\in N_{n}\}$ is a $N_{n}$ regular conditional distribution of $X_{n+1}$ given $\mathbf{X}$. Assume also that $\{\mu_{n+1}^{(\mathbf{x},x_{n+1})},\ (\mathbf{x},x_{n+1})\in N_{n}\times \R_{+}\}$ is a $N_{n+1}:=N_{n}\times \R_{+}$ regular conditional distribution of $X_{n+2}$ given $(\mathbf{X},X_{n+1})$.
 Then we have:
\begin{description}
\item {a) } A $N_{n}$ regular conditional distribution of $(X_{n+1},X_{n+2})$ given $\mathbf{X}$ is given by the composition kernel $\nu:=\mu_{n+1}\otimes\mu_{n+2}$ defined in (\ref{condos}).

\item {b) } Take a family of bidimensional random variables $\{(Z_{n+1}^{\mathbf{x}},Z_{n+2}^{\mathbf{x}}),\  \mathbf{x}\in N_{n}\}$, such that for each $\mathbf{x}\in N_{n}$,  the distribution of $(Z_{n+1}^{ \mathbf{x}}, Z_{n+2}^{ \mathbf{x}})$ is given by
$\nu_{n}^{\mathbf{x}}$, as in part a).  For each $\mathbf{x}\in N_{n}$, the family  $\{\mu_{n+1}^{(\mathbf{x},x_{n+1})},\ x_{n+1}\in \R_{+}\}$ is a $ \R_{+}$ regular conditional distribution of  $Z_{n+2}^{\mathbf{x}}$ given $Z_{n+1}^{\mathbf{x}}$. \end{description}\label{prvers}\end{prop}

 Now we are in a position to prove one of the main results of this section.
 \begin{prop}
Let $  \{ N(t): t \geq 0 \} $ be a general counting process with inter renewal epochs $\seq{X}$. Let $T$ be a random time independent from the process. Assume that:
\begin{description}
\item {a) }$T$ is a DFR random variable.
\item {b) }$T$ and $X_{1}$ don't have simultaneously positive mass at 0.
\item {c)} For each $n=1,2,3,\dots,$ there exists a $ N_{n}$ distributional version of $X_{n+1}$ given $(X_{1},\dots,X_{n})$ (say $\{Z_{n+1}^{\mathbf{x}},\quad  \mathbf{x}\in N_{n}\}$) verifying that
   \begin{description}\item {c.1) }$Z_{2}^{x_{1}}\leq_{ST}X_{1},\quad x_{1}\in  N_{1},$
\item{c.2) } For $n=1,2,\dots,$ $ N_{n+1}\subseteq N_{n}\times \R^{+}$ and
\[Z_{n+2}^{(\mathbf{x},x_{n+1})}\leq_{ST} Z_{n+1}^{\mathbf{x}},\quad \hbox{for all} \quad (\mathbf{x},x_{n+1})\in N_{n+1}.\]
    \end{description}

\end{description}
Then $N(T)$ is discrete DFR. \label{t2*}

\end{prop}

{\bf Proof:}  The result will follow by checking (\ref{dfrn0}) and (\ref{dfrn}).
To the first inequality, note that assumption c.1), allow us to apply Proposition \ref{prcon},  thus having
\[ E^{2} (\overline {F}_T (X_{1}))\leq E (\overline {F}_T (X_{1}+X_{2}))\overline {F}_T(0)\leq E (\overline {F}_T (X_{1}+X_{2})),\]
and therefore  (\ref{dfrn0}) holds true. As for (\ref{dfrn}), thanks to c.2), we have, for fixed $n$, the family of random variables  $\{Z_{n+2}^{(\mathbf{x},x_{n+1})},\  (\mathbf{x},x_{n+1})\in N_{n+1}\}$, being a $N_{n+1}$ distributional version of $X_{n+2}$ given $(X_{1},\dots,X_{n+1})$.  We extend it to a $N_{n}\times \R^{+}$ version by taking $Z_{n+2}^{(\mathbf{x},x_{n+1})}=0$ for each $(\mathbf{x},x_{n+1})\in N_{n}\times \R^{+}$ such that $(\mathbf{x},x_{n+1})\not \in N_{n+1}$. Then, we have
\begin{equation}Z_{n+2}^{(\mathbf{x},x_{n+1})}\leq_{ST} Z_{n+1}^{\mathbf{x}} \hbox{ for all }(\mathbf{x},x_{n+1})\in N_{n}\times \R^{+}.\label{monoto}\end{equation}
This version, together with $\{Z_{n+1}^{\mathbf{x}},\  \mathbf{x}\in N_{n}\}$, the $N_{n}$ distributional version of $X_{n+1}$ given $(X_{1},\dots,X_{n})$, provides us, applying Proposition \ref{prvers} a), a $N_{n}$ regular conditional distribution of $(X_{n+1},X_{n+2})$ given $\mathbf{X}$, with its corresponding distributional version  $\{(Z_{n+1}^{\mathbf{x}},Z_{n+2}^{\mathbf{x}}),\  \mathbf{x}\in N_{n}\}$.

By Proposition  \ref{prvers} b) $\{Z_{n+2}^{(\mathbf{x},x_{n+1})}, \quad \ x_{n+1}\in \R^{+}\}$ is a  $\R^{+}$ version of $Z_{n+2}^{\mathbf{x}}$ given $Z_{n+1}^{\mathbf{x}}$, which verifies (\ref{monoto}), so that we can apply Proposition \ref{prcon} to the random vector $(X_{1},X_{2})=:(Z_{n+1}^{\mathbf{x}},Z_{n+2}^{\mathbf{x}})$ and any arbitrary DFR survival function $\overline {H}$ such that $\overline {H}(0)=1$, thus having,
\begin{equation}E^{2}[ \overline {H} (Z_{n+1}^{\mathbf{x}})] \leq E[ \overline {H}( Z_{n+1}^\mathbf{x}+Z_{n+2}^{\mathbf{x}})].\label{conass2}\end{equation}
Thus, conditions in Proposition \ref{prdfr2} are verified, and (\ref{dfrn}) holds true.  This completes the proof of Proposition \ref{t2*}.$\quad \Box$

\section{Applications to maintenance models}
Maintenace models describe the behaviour of a mechanism subject to repair.  When the mechanism fails a repair is performed.  This repair can restore the system to its initial state (perfect maintenace) or to the condition as it was before the failure occurred (minimal maintenace). Also, it can leave the system in an intermediate state between both of them (imperfect maintenace), or even in an operational state worse than just before the failure (worse maintenace). This kind of models were introduced by Kijima and Sumita (see \cite{kisuau}). For a recent review of these models we refer to \cite{inma}. Among them, the so-called virtual age models, associate, to each time t, the virtual age of the system $v(t)$, which represents the operational stage of the system after successive failures and its corresponding repairs, in such a way that, if $\bar{G}$ represents the lifetime of a new unit, then the lifetime of a unit with virtual age $v$ is given by
\begin{equation}\bar{G}(z\vert v):=\left\{
                                           \begin{array}{ll}
                                             \frac{\overline{G}(v+z)}{\overline{G}(v)}, & \hbox{if $\overline{G}(v)>0$ ;} \\
                                             0, & \hbox{if $\overline{G}(v)=0$.}
                                           \end{array}
                                         \right.,\quad v\geq0.
\label{survir}\end{equation}
In virtual age models the underlying renewal process is defined through the interarrival times $\seq{X}$ representing the time elapsed between consecutive failures of the mechanism. We also have the sequence $\seq{V}$ representing the virtual age of the system just after the n-th repair (we assume there is no delay due to the repair).  The virtual age of the system usually depends on $\seq{X}$ and a family of nonnegative random variables $\seq{A}$, representing the degree of repair performed at $n$. The value $A_{n}=0$ represents the best possible repair of the model under consideration, whereas $A_{n}=1$ stands for minimal maintenace (the repair leaves the system as it was just before the failure). For instance, Kijima (cf. \cite{kisome}) considered two models, in which the virtual age is defined as  follows:
\begin{eqnarray}
&V_{n}=V_{n-1}+A_{n}\cdotp X_{n},\quad &\hbox{Kijima type I},\label{kijit1}\\
&V_{n}=A_{n}\cdotp(V_{n-1}+ X_{n}),\quad &\hbox{Kijima type II}.\label{kijit2}
\end{eqnarray}
In general, the virtual age is defined, for $n=1,2,\dots $, as
\begin{equation}V_{n}:=\Psi(V_{n-1},X_{n},A_{n})=\Psi_{n}(X_{1},\dots,X_{n},A_{1},\dots,A_{n}), \quad V_{0}=0,\label{virage}\end{equation}
in which $\Psi$ is a suitable function.  Once the virtual age at instant $n$ is known, the survival function of the next repair is evaluated as follows.
\begin{equation} P(X_{n+1}>z|V_{n}=v)=\bar{G}(z\vert v),\quad v\geq 0.\label{convva}\end{equation}
In the particular case of  $\seq{A}$ being deterministic nonnegative numbers, it is clear by (\ref{virage}) and (\ref{convva}) that a $R_{+}^{n}$ distributional version of $X_{n+1}$ given $X_{1},\dots,X_{n}$ is given by a family of random variables $\{Z_{n+1}^{\mathbf{x}},\quad  \mathbf{x}\in R_{+}^{n}\}$, whose reliability function is given by
\begin{equation} P(Z_{n+1}^{\mathbf{x}}>z)=\bar{G}(z\vert \Psi_{n}(\mathbf{x},A_{1},\dots,A_{n})),\quad z\geq 0.\label{convv2}\end{equation}
Of particular interest is the case in which $A_{n}=q, \ n=1,2,\dots$, for a given $0\leq q\leq 1$.  Then, we obtain the generalized renewal processes considered by Kijima and Sumita in (\cite{kisuau}).

In the following Proposition we prove the d-DFR property for $N(T)$, for general Kijima type I models. The sequence $\seq{A}$ can be random, although assume that, for each $n$, $A_{n+1}$ is independent of $X_{1},\dots,X_{n}$.  This assumption is also made in the original paper of Kijima  \cite{kisome}, along with the assumption that the sequence $\seq{A}$ is formed of independent random variables, as well as $0\leq A_{n}\leq 1$.  We don't need to make these last assumptions, as the method of proof is practically the same.
This result can be useful when an age replacement is used up to a random time $T$, in which a preventive maintenance can be performed.

\begin{prop}
Let $  \{ N(t): t \geq 0 \} $ be the counting process associated to a Kijima type I repair model as in (\ref{kijit1}). Assume that the sequence of random repairs $\seq{A}$ satisfies that $A_{n+1}$ is independent of $(X_{1},\dots,X_{n}), \ n=1,2,\dots$. Let $T$ be a random time independent from the process. If the following conditions hold true:
\begin{description}
\item {a) } $T$ is a DFR random variable,
\item {b) } $X_{1}$ is increasing failure rate.
\end{description}
Then, $N(T)$ is discrete DFR.
\label {prgenr}
\end{prop}
{\bf Proof:} First of all, we will consider that the values in sequence $\seq{A}$ is formed of deterministic numbers. The result will follow applying Proposition \ref{t2*}.  As $X_{1}$ is IFR, then it doesn't have positive mass at 0 (see (\cite[p.\ 109]{maolli})), so that condition b) is satisfied.  As for c), note that, taking into account (\ref{kijit1}), in this case $\psi_{n}$ in (\ref{virage}) has the form
\begin{equation}\Psi_{n}(X_{1},\dots,X_{n},A_{1},\dots,A_{n})=A_{1}X_{1}+\dots+ A_{n}X_{n}.\label{ps1ki1}\end{equation}
As $X_{1}$ is IFR, then $\bar{G}(\cdotp\vert y)$, as in (\ref{survir}) is decreasing in $y$, and therefore, taking into account (\ref{convva}) and (\ref{ps1ki1}), we have
\begin{equation}P(Z_{1}^{x_{1}}>z)=\bar{G}(z\vert A_{1}x_{1})\leq \bar{G}(z\vert 0) =P(X_{1}>z),\label{firsin}\end{equation}
 so that $Z_{1}^{x}\leq_{ST} X_{1}$, and c.1) is verified.  c.2) is checked in a similar way.  In this case we have
\begin{eqnarray*}P(Z_{n+2}^{\mathbf{x},x_{n+1}}>z)&=&\bar{G}(z\vert A_{1}x_{1}+\dots +A_{n+1}x_{n+1})]\\
&\leq& \bar{G}(z\vert A_{1}x_{1}+\dots +A_{n}x_{n})]=P(Z_{n+1}^{\mathbf{x}}>z)\end{eqnarray*}
and $Z_{n+2}^{\mathbf{x},x_{n+1}}\leq_{ST}Z_{n+1}^{\mathbf{x}}$.  Thus, conditions in Proposition \ref{t2*} are verified when the sequence $\seq{A}$ is formed of a sequence of constant values and in this case, $N(T)$ is d-DFR.  Now, assume that the sequence $\seq{A}$ is formed of random variables.  Fix $n$ a natural number.  In order to show (\ref{disdfr}) we will condition to the random vector $\mathbf{\Theta}_{n}:=(A_{1},\dots,A_{n+1})$, that is
\begin{equation}\kern-10pt P(N(T)\geq n+i)=\int_{\R^{n+1}_{+}}P(N(T)\geq n+i|\mathbf{\Theta}_{n}=\mathbf{\theta})dF_{\mathbf{\Theta_{n}}}(\theta),\  i=0,1,2.\label{mixtur}\end{equation}
For fixed $\theta:=(a_{1},\dots,a_{n+1})$, consider a Kijima type I model $  \{ N_{\theta}(t): t \geq 0 \} $ as in (\ref{kijit1}), in which $\seq{A}$ verifies $A_{i}=a_{i},\ i=1,2,\dots,n+1$ (the rest of values in the sequence can be defined arbitrarily).  Note that, taking into account (\ref{nocomm}), $P(N(T)\geq n)$ depends only on $(X_{1},\dots,X_{n})$, whose distribution is specified in terms of $X_{1}$ and the virtual ages $V_{1}=A_{1}X_{1},\dots,V_{n-1}=A_{1}X_{1}+\dots+A_{n-1}X_{n-1})$, so that using the fact that each $A_{i}$ is independent of the process of past lifetimes $X_{1},\dots,X_{n-1}$,
\begin{equation}P(N(T)\geq n+i|\mathbf{\Theta}_{n}=\theta)=P(N_{\theta}(T)=n+i),\quad i=0,1,2,\label{identi}\end{equation}
and (\ref{mixtur}) and (\ref{identi}) gives us
\begin{equation}P(N(T)\geq n+1)=\int_{\R^{n+1}_{+}}P(N_{\theta}(T)\geq n+1) dF_{\mathbf{\Theta_{n}}}(\theta),\label{casica}\end{equation}
By the previous equality and taking into account that $  \{ N_{\theta}(t): t \geq 0 \} $ is a Kijima type I model with deterministic repairs (and therefore, verifying the d-DFR property) we obtain
\begin{eqnarray*}&&P(N(T)\geq n+1)\leq \int_{\R^{n+1}_{+}}P^{\frac{1}{2}}(N_{\theta}(T)\geq n+2)P^{\frac{1}{2}}(N_{\theta}(T)\geq n)dF_{\mathbf{\Theta_{n}}}(\theta)\\&&\leq\left(\int_{\R^{n+1}_{+}}P(N_{\theta}(T)\geq n+2)dF_{\mathbf{\Theta_{n}}}(\theta)\right)^{\frac{1}{2}}\left(\int_{\R^{n+1}_{+}}P(N_{\theta}(T)\geq n)dF_{\mathbf{\Theta_{n}}}(\theta)\right)^{\frac{1}{2}}\end{eqnarray*}
Now, we use the previous inequality together with (\ref{mixtur}) and (\ref{identi}) to conclude
\[P(N(T)\geq n+1)\leq P^{\frac{1}{2}}(N(T)\geq n+2)P^{\frac{1}{2}}(N(T)\geq n)\quad \Box\]

\begin{rem}
The previous result is more general than the one given in \cite[Corollary 5.1.]{basath}, in two senses.  Firstly, we allow the repairing sequence to be random and secondly, the initial repair $X_{1}$ is more general. The question is that in \cite{basath} we used an analogous proof than the one considered in \cite{german} to obtain IFR preservation results, and for this property stronger conditions on the arrival epoch are needed.

Also, it is interesting to point out that, under the same conditions as in Proposition \ref{prgenr}, Kijima type II models don't verify the DFR property for $N(T)$.  To show that, consider a Kijima type II model, in which $X_{1}$ is increasing failure rate, and $A_{1}=1$, $A_{2}=0$ (the problem here is that the second interarrival time decreases with respect the first one, but, as the second repair is perfect, the process at this point 'restarts as new').  Consider that $T$ is an exponential random variable of mean 1 (that is, $ \overline {F}_T(x)=e^{-x},\ x\geq 0$). Firstly, we observe in (\ref{kijit1}), (\ref{kijit2}) that the virtual age $V_{1}$ in Kijima type I and Kijima type II coincide (as $V_{0}=0$) and, taking into account (\ref{firsin}) and Proposition \ref{prcon}, the discrete DFR property for $n=1$ is satisfied, that is,
\begin{equation} E^2[ e^{-X_{1}}] \leq E[ e^{- (X_{1}+X_{2})}]. \label{counte}\end{equation}
On the other hand, as $A_{2}=0$, we see from (\ref{kijit2}) that $V_{2}=0$, so that $X_{3}=_{d}X_{1}$, and is independent of $X_{1}$ and $X_{2}$.  Thus, we can write
\[ E[ e^{- S_{3}}]=E[ e^{- (X_{1}+X_{2})}]E[e^{-X_{3}}]=E[ e^{- (X_{1}+X_{2})}]E[e^{-X_{1}}].\]
Assume that the discrete DFR property is true for $n=2$.  Taking into account (\ref{dfrn}) and the previous expression we have
\[E^{2}[ e^{- (X_{1}+X_{2})}]\leq E[ e^{- S_{3}}]E[e^{-X_{1}}]=E[ e^{- (X_{1}+X_{2})}]E^{2}[e^{-X_{1}}].\]
We conclude therefore that,
\begin{equation}E[ e^{- (X_{1}+X_{2})}]\leq E^{2}[e^{-X_{1}}].\label{count2}\end{equation}
Thus, taking into account (\ref{counte}) and (\ref{count2}), the d-DFR for $n=1$ and $n=2$ would imply
\[E[ e^{- (X_{1}+X_{2})}]=E^{2}[e^{-X_{1}}],\]
which is not, in general, true.  Take for instance $X_{1}$ to be uniformly distributed on $(0,1)$.  Then, we can check easily that
\[ P(Z_{2}^{x_{1}}>z)=\bar{G}(z\vert x_{1})=\frac{1-x_{1}-z}{1-x_{1}},\quad 0\leq z\leq 1-x_{1}\]
and thus, the distribution of $Z_{2}^{x_{1}}$ is uniform on $(0,1-x_{1})$.  Then, we can write
\[E[ e^{-(X_{1}+X_{2})}]=\int_{0}^{1}dx_{1} e^{-x_{1}}\int_{0}^{1-x_{1}}dx_{2}\frac{e^{-x_{2}}}{1-x_{1}}=\int_{0}^{1}dx_{1} e^{-x_{1}}\int_{0}^{1}du e^{-(1-x_{1})u},\]
which is clearly greater than
\[E^{2}[ e^{-X_{1}}]=\int_{0}^{1}dx_{1} e^{-x_{1}}\int_{0}^{1}du e^{-u}.\]

\end{rem}

\section{The DFR property under association of the interarrival times}
In  \cite{basath}, Proposition 3.1 (a) sufficient conditions for two random variables in order to satisfy (\ref{dfrn0}) were given.  These conditions involve the concept of association (cf \cite{esprwa}).  We now recall this concept.

\begin{defn} A random vector $\mathbf{X}:=(X_{1},\dots,X_{n})$ is said to be \emph{associated} if
\begin{equation}Cov(f(\mathbf{X}),g(\mathbf{X}))=E[f(\mathbf{X})\cdotp g(\mathbf{X})]-E[f(\mathbf{X})]\cdotp E[g(\mathbf{X})]\geq0\label{poscor}\end{equation}
for all increasing functions $f,g:\R^{n}\rightarrow \R$ such that the expectation exists.
\end{defn}
On the other hand, we know from \cite{basath} that
\begin{prop}[ \cite{basath}, Proposition 3.1(a)] Let $\seq{X}$ be an arbitrary sequence of interarrival times and $T$ a random time independent of them.  Assume that $T$ is a DFR random variable. If $(X_{1},X_{2})$ is associated and $X_{2}\leq_{ST}X_{1}$, then
\[E^2[ \overline {F}_T (X_{1})] \leq E[ \overline {F}_T (S_{2})]\overline {F}_T (0) . \]
\label{prasso}
\end{prop}

As an immediate consequence of the previous result and Proposition \ref{prdfr2} we have the following

\begin{cor}Let $  \{ N(t): t \geq 0 \} $ be a general counting process with inter renewal epochs $\seq{X}$. Let $T$ be a random time independent from the process. Assume that
\begin{description}
\item {a) }$T$ is a DFR random variable
\item {b) }$T$ and $X_1$ don't have 0 as a common discontinuity point.
\item {c) }$(X_{1},X_{2})$ is associated and $X_{2}\leq_{ST}X_{1}$.
\item {d) }  For fixed $n$ there exists an $ N_{n}$ distributional version of $(X_{n+1},X_{n+2})$ given $(X_{1},\dots,X_{n})$ (say $\{(Z_{n+1}^{\mathbf{x}},Z_{n+2}^{\mathbf{x}}),\quad  \mathbf{x}\in N_{n}\}$ satisfying:\[(Z_{n+1}^{\mathbf{x}},Z_{n+2}^{\mathbf{x}})\hbox{ is associated and } Z_{n+2}^{\mathbf{x}}\leq_{ST}Z_{n+1}^{\mathbf{x}},\quad \hbox{for all } \mathbf{x}\in N_{n}.\]
    \end{description}
    Then, $N(T)$ is discrete DFR. \label{cassoc}
    \end{cor}

    {\bf Proof:}  We will show that, under these hypothesis, (\ref{dfrn0}) and  (\ref{dfrn}) are verified.    In fact, Condition c)  and Proposition \ref{prasso} imply immediately (\ref{dfrn0}).  To prove (\ref{dfrn}), fix $n$ and consider $\bar{H}$ a DFR survival function with $\bar{H}(0)=1$.  Condition d), together with Proposition \ref{prasso} applied to $(X_{1},X_{2}):=(Z_{n+1}^{\mathbf{x}},Z_{n+2}^{\mathbf{x}})$ implies that
    \[E^2[ \overline {H} (Z_{n+1}^{\mathbf{x}})] \leq E[ \overline {H} (Z_{n+1}^{\mathbf{x}}+Z_{n+2}^{\mathbf{x}})]. \]
    This, together with Proposition \ref{prdfr2}, proves (\ref{dfrn}), and therefore the d-DFR property for $N(T)$.$\quad \Box$

    \begin{rem} In \cite{basath}, Proposition 3.2.  the same condition (c) was required, but it was also required the association of the random vector $(S_{n},X_{n+1},X_{n+2})$, for $n=1,2,\dots$.  It is not hard to find a renewal process satisfying the conditions in Corollary \ref{cassoc}, but not verifying the previous association condition. For instance, take $Y$ a random variable with support on $[1, \infty )$,  and let $\seq{W}$ be a sequence of independent and identically distributed non-negative random variables independent of $Y$ and verifying $Y \geq_{\hbox {ST}} W_{1}$. Define the inter-renewal epochs as follows
\begin{equation}
X_1 =Y \quad X_2=W_{1} \quad \hbox{and} \quad
X_{n} =\frac{W_{n-1}}{Y},\ n=3,4,\dots \label{coun21}\end{equation}
It is easy to check that conditions in Corollary \ref{cassoc} are verified.  However, $(X_{1},X_{2},X_{3})$ is not associated.  To this end, take $Y$ verifying that
\begin{equation}E[Y]\cdotp E\left[\frac{1}{Y}\right]>1\label{coun22}\end{equation}
(this inequality is satisfied, for instance, if $Y=1+Ber(p)$ with $Ber(p)$ a Bernoulli random variable with probability of success $0<p<1$).
Now, consider that $Y$ and $W_{1}$ have finite mean. We will show that $(X_{1},X_{3})$ is not associated.  To this end, take the increasing functions \[f(x,y)=x,\quad g(x,y)=y, \qquad (x,y)\in \R^{2}.\] Taking into account (\ref{coun21}),  (\ref{coun22}), the fact that $Y$ and $W_{2}$ are independent and that $E[W_{2}]>0$, we have
\begin{eqnarray*}E[f(X_{1},X_{3})g(X_{1},X_{3})]&=&E[W_{2}]<E[Y]E[W_{2}]E\left[\frac{1}{Y}\right]\\&=&E[Y]E\left[\frac{W_{2}}{Y}\right]=E[f(X_{1},X_{3})]\cdotp E[g(X_{1},X_{3})]\end{eqnarray*}
Thus $(X_{1},X_{3})$ is not associated and therefore $(X_{1},X_{2},X_{3})$ is not associated (see  \cite[p.123]{mustco}).
 In general, the association of random variables does not imply the conditional association and vice versa (see \cite{yuyaco} for some examples)) so that the conditions given in Corollary \ref{cassoc} and in \cite{basath}, Proposition 3.2. are different.  The practical advantage in this last result is that we don't need to build the conditional distributions of the interarrival times.\end{rem}

\appendix
\section{Appendix}
In this Appendix we include the proof of Proposition \ref{prvers}

\medskip

{\bf Proof of Proposition \ref{prvers}} Let $\{\nu_{n}^{\mathbf{x}},\ \mathbf{x}\in N_{n}\}$ be the composition kernel, as defined in part (a).  To show that this kernel is a $N_{n}$ regular conditional distribution of $(X_{n+1},X_{n+2})$ given $(X_{1},\dots,X_{n})$,  we need to show for each  $A\in {\mathcal B}(N_{n})$ and $B\in {\mathcal B}(\R_{+}^2)$(recall (\ref{promed}))
\begin{equation}\int_{A\times B } dF_{n+2}(\mathbf{x},x_{n+1},x_{n+2})=\int_{A}dF_{n}(\mathbf{x})\int_{B} d\nu_{n
}^{\mathbf{x}}(x_{n+1},x_{n+2}).\label{basic}\end{equation}
Firstly, we will verify the previous equality for Borel sets $A\in N_{n}$,  and $B$ of the form
\begin{equation}B:=B_{1}\times B_{2},\quad B_{i}\in {\mathcal B}(\R_{+}),\quad i=1,2\label{rectan}\end{equation}
and by a classical extension argument we will conclude that (\ref{basic}) is also verified for Borel sets $B\in \R^{2}_{+}$.
To show (\ref{basic}), note firstly that, by assumption, $\{\mu_{n+1}^{(\mathbf{x},x_{n+1})},\ (\mathbf{x},x_{n+1})\in N_{n}\times \R_{+}\}$ correspond to a  $N_{n}\times \R_{+}$ regular conditional distribution of $X_{n+2}$ given $(\mathbf{X},X_{n+1})$, so that we can apply (\ref{promed}) to this version, and write for $B$ as in (\ref{rectan})
\begin{eqnarray}&&\int_{A\times B_{1}\times B_{2}}dF_{n+2}(\mathbf{x},x_{n+1},x_{n+2})\nonumber\\&&=\int_{A\times B_{1}}dF_{n+1}(\mathbf{x},x_{n+1})\int_{B_{2}}d\mu_{n+1}^{(\mathbf{x},x_{n+1})}(x_{n+2})
.\label{basi1}\end{eqnarray}
Consider now the function
\begin{eqnarray*}g(\mathbf{x},x_{n+1})&=&\int_{B_{2}}d\mu_{n+1}^{(\mathbf{x},x_{n+1})}(x_{n+2})\cdotp 1_{A\times B_{1}}(\mathbf{x},x_{n+1})\\&=&\int_{B_{2}}d\mu_{n+1}^{(\mathbf{x},x_{n+1})}(x_{n+2})\cdotp
1_{A}(\mathbf{x})\cdotp 1_{B_{1}}(x_{n+1}),\end{eqnarray*}
As by hypothesis $\{\mu_{n}^{\mathbf{x}},\ \mathbf{x}\in N_{n}\}$ correspond to a $ N_{n}$ distributional version of $X_{n+1}$ give $\mathbf{X}$, we can  apply (\ref{prodint}) to the previous function and write the right-hand side term in (\ref{basi1}) as
\begin{eqnarray} &&\int_{A\times B_{1}}dF_{n+1}(\mathbf{x},x_{n+1})\int_{B_{2}}d\mu_{n+1}^{(\mathbf{x},x_{n+1})}(x_{n+2})=Eg(\mathbf{X},X_{n+1})\nonumber\\
&&=\int_{N_{n}}dF_{n}(\mathbf{x})\int_{\R_{+}}d\mu_{n}^{\mathbf{x}}(x_{n+1}) g(\mathbf{x},x_{n+1})
\nonumber\\&&=\int_{A}dF_{n}(\mathbf{x})\int_{B_{1}}d\mu_{n}^{\mathbf{x}}(x_{n+1}) \int_{B_{2}}d\mu_{n+1}^{(\mathbf{x},x_{n+1})}(x_{n+2})\nonumber\\
&&=\int_{A}dF_{n}(\mathbf{x})\int_{B_{1}\times B_{2}}d\nu_{n}^{\mathbf{x}}(x_{n+1},x_{n+2}),\label{basi2}\end{eqnarray}
where, in the last equality we have used (\ref{condos}).  Therefore, (\ref{basi1}) and (\ref{basi2}) show (\ref{basic}) for a fixed Borel set $A\in N_{n}$ and sets of the form (\ref{rectan}).  Using a classical monotone class argument (cf. \cite[Thm. 1.1., p.2]{kafoun})
(\ref{basic}) holds true for all Borel sets $B\in \R_{+}^{2}$.  This means that the family $\{\nu_{n}^{\mathbf{x}},\ \mathbf{x}\in N_{n}\}$ satisfies Definition \ref{defcon}, thus being a $N_{n}$ regular conditional distribution of $(X_{n+1},X_{n+2})$ given $\mathbf{X}$, and part (a) is completed. Part (b) follows easily by (a).  Let $\{(Z_{n+1}^{\mathbf{x}},Z_{n+2}^{\mathbf{x}}),\ \mathbf{x}\in N_{n}\}$ be the family of random variables as  given in (b). First of all, note that for fixed $\mathbf{x}\in N_{n}$, $\{\mu_{n+1}^{(\mathbf{x},x_{n+1})},\ x_{n+1}\in R_{+}\}$ is a probability kernel from $R_{+}$ to $R_{+}$ . Secondly note that (\ref{condos}) gives (\ref{promed}) for $(Z_{n+1}^{\mathbf{x}},Z_{n+2}^{\mathbf{x}})$ in terms of the marginal of the first variable ($Z_{n+1}^{\mathbf{x}}$) and $\{\mu_{n+1}^{(\mathbf{x},x_{n+1})},\ x_{n+1}\in R_{+}\}$, so that Definition 5 is verified and we conclude that the family  $\{\mu_{n+1}^{(\mathbf{x},x_{n+1})},\ x_{n+1}\in \R_{+}\}$ is a $ \R_{+}$ regular conditional distribution of  $Z_{n+2}^{\mathbf{x}}$ given $Z_{n+1}^{\mathbf{x}}$, as claimed.

\section*{Aknowledgements}

This work has been supported by
research projects MTM2010-15311 and by FEDER
funds. The first and second authors acknowledge the support of DGA S11 and E64,
respectively.

\end{document}